\documentclass[11pt]{article}

\usepackage{amsmath,amsfonts,amsthm,amssymb,graphicx,tikz,tikz-cd,url,hyperref,enumerate}
\usepackage[all,2cell,ps]{xy}
\usepackage{dsfont}
\usepackage[symbol]{footmisc}
\usepackage{overpic}
\usepackage{thm-restate}
\theoremstyle{plain}
\newtheorem{thm}{Theorem}[section]

\newtheorem{lem}[thm]{Lemma}
\newtheorem{prop}[thm]{Proposition}

\newtheorem{cor}[thm]{Corollary}

\newtheorem{ques}[thm]{Question}

\theoremstyle{definition}

\newtheorem{rem}[thm]{Remark}

\theoremstyle{remark}

\newcommand{\C}{\mathbb{C}}

\newcommand{\F}{\mathbb{F}}

\newcommand{\Hy}{\mathbb{H}}

\newcommand{\N}{\mathbb{N}}

\newcommand{\R}{\mathbb{R}}

\newcommand{\Z}{\mathbb{Z}}

\newcommand{\sig}{\sigma}

\DeclareMathOperator{\SL}{SL}
\DeclareMathOperator{\PSL}{PSL}

\DeclareMathOperator{\PGL}{PGL}

\DeclareMathOperator{\SO}{SO}
\DeclareMathOperator{\MSO}{MSO}

\DeclareMathOperator{\id}{id}

\DeclareMathOperator{\Isom}{Isom}
\DeclareMathOperator{\Aut}{Aut}
\DeclareMathOperator{\Inn}{Inn}
\DeclareMathOperator{\Out}{Out}
\DeclareMathOperator{\Hom}{Hom}

\DeclareMathOperator{\Mod}{Mod}
\DeclareMathOperator{\Diff}{Diff}

\DeclareMathOperator{\ord}{ord}

\usepackage{tikz-cd}
 \allowdisplaybreaks

\usepackage{color}

\definecolor{orange}{rgb}{1,0.5,0}

\hypersetup{
  colorlinks,
  citecolor=black,
  linkcolor=black,
  urlcolor=black}

\usetikzlibrary{arrows}

\title{On signatures of the atoroidal bundles of Kent--Leininger}
\author{Jean-Fran\c cois Lafont, Nicholas Miller, and Lorenzo Ruffoni}

\date{}

\newcommand{\Addresses}{{% additional braces for segregating \footnotesize
\bigskip
\footnotesize

\textsc{Department of Mathematics, Ohio State University, Columbus, OH 43210}\par\nopagebreak
\textit{E-mail address}, J.-F. Lafont: \texttt{jlafont@math.ohio-state.edu} \\

\textsc{Department of Mathematics, University of Oklahoma, Norman, OK 73019}\par\nopagebreak
\textit{E-mail address}, N. Miller: \texttt{nickmbmiller@ou.edu} \\

\textsc{Department of Mathematics, Binghamton University, Binghamton, NY 13902}\par\nopagebreak
\textit{E-mail address}, L. Ruffoni: \texttt{lruffoni@binghamton.edu} 
}}
\begin{document}

\maketitle

%%%%%%%%%%%%%%%%%%%%
%%%%%%%%%%%%%%%%%%%%
%%%%%%%%%%%%%%%%%%%%

\begin{abstract}
\noindent We show that there are infinitely many homeomorphism types of atoroidal surface bundles over surfaces which have signature zero.
\end{abstract}

%%%%%%%%%%%%%%%%%%%%
%%%%%%%%%%%%%%%%%%%%
%%%%%%%%%%%%%%%%%%%%

\section{Introduction}
In $3$-manifold topology, it is well known that one can construct surface bundles over the circle by taking the mapping torus of a homeomorphism $\phi$ of a closed or punctured surface.
When the mapping class associated to $\phi$ is pseudo-Anosov, seminal work of Thurston \cite{Otal,Otal2} shows that this mapping torus is moreover a finite volume, hyperbolic manifold.
In fact, Thurston shows the stronger fact that the hyperbolicity of such a mapping torus is equivalent to it being atoroidal, that is, not containing an essential torus.

Motivated by the desire to extend this picture to dimension four, for several decades it has been an open question whether there similarly exist atoroidal surface bundles over surfaces in dimension four.
Such manifolds represent natural candidates for constructing hyperbolic manifolds in the spirit of Thurston.
This year, groundbreaking work of Kent--Leininger \cite{KL} resolves this by showing that a plethora of atoroidal surface bundles over surfaces exist.
In what follows, $S_g$ always denotes a closed surface of genus $g\ge 2$, and $\Mod(S_g)$ denotes its mapping class group.

%%%%%%%%%%%%%%%%%%%%
%%%%%%%%%%%%%%%%%%%%

\begin{thm}[Kent--Leininger]\label{thm:KL}
For every $g\ge 4$, there exist infinitely many distinct commensurability classes of purely pseudo-Anosov subgroups $\pi_1(S_h)<\Mod(S_g)$ where $h\ge 2$. 
In particular, there are infinitely many homeomorphism types of atoroidal surface bundles over surfaces.
\end{thm}

%%%%%%%%%%%%%%%%%%%%
%%%%%%%%%%%%%%%%%%%%

Recall that a purely pseudo-Anosov subgroup of the mapping class group is one for which every non-trivial element is a pseudo-Anosov mapping class. 
The equivalence of the existence of such subgroups and the existence of an atoroidal surface bundle over a surface is well known; see Section \ref{sec:surfacebundle} for details.
Though not stated explicitly in \cite{KL}, the statement on homeomorphism types follows from the fact that a bundle with base $S_h$ and fiber $S_g$ has Euler characteristic $4(g-1)(h-1)$, and both $g$ and $h$ can be chosen to be arbitrarily large.
Furthermore, it is known that surface bundles over surfaces have finitely many fibering structures \cite{Johnson}, in contrast to the $3$-dimensional setting \cite{Thurston}.

The work of Kent--Leininger opens the door to exploring whether such bundles admit interesting geometric structures. 
For instance, the following are well-known and important open questions that have been widely circulated over the past few decades.

%%%%%%%%%%%%%%%%%%%%
%%%%%%%%%%%%%%%%%%%%

\begin{ques}\label{ques:main}
Suppose that $M$ is an atoroidal surface bundle over a surface.
Is $\pi_1(M)$ Gromov hyperbolic?
Does $M$ possess a negatively curved Riemannian metric?
Does $M$ possess a hyperbolic metric?
\end{ques}

%%%%%%%%%%%%%%%%%%%%
%%%%%%%%%%%%%%%%%%%%

In this paper, by hyperbolic, we mean real hyperbolic. Note that this is the only locally symmetric structure that such $M$ can admit.
For instance, if a closed manifold is modeled on $\Hy^2\times \Hy^2$, then its fundamental group must contain $\Z^2$ by \cite[Cor 2.9]{PR72}.
On the other hand, surface bundles are known not to admit complex hyperbolic structures by \cite{Kapovich}. 
Alternatively, for the specific kind of surface bundles considered in this paper,  one can argue that if $M$ is a closed $4$-manifold with a complex hyperbolic metric, then $\chi(M)=3|\sig (M)|$ by \cite{VI85}, but the Euler characteristic of an $S_g$-bundle over $S_h$ is $4(g-1)(h-1)$, while we compute that the signature vanishes.
% must be odd, because  $H^2(M;\R)\cong\R$ by Berger \cite[\S 3]{BE60}. But the Euler characteristic of an $S_g$-bundle over $S_h$ is even, as remarked above. See also  for an alternative proof.

We also remark that in Question \ref{ques:main}, an affirmative answer to any of the questions listed implies an affirmative answer to their predecessors, therefore $M$ being a real hyperbolic manifold is the strongest statement one could ask for.
To our knowledge, it is widely expected that such $M$ should {\bf not} be real hyperbolic manifolds, although $\pi_1(M)$ is expected to be Gromov hyperbolic.
In the special case of the manifolds arising from the construction of \cite{KL}, Kent--Leininger do indeed conjecture the latter is always the case.

Producing or, alternatively, ruling out hyperbolic structures on closed manifolds $M$ can be difficult in the absence of obvious obstructions to hyperbolicity (e.g., being aspherical and atoroidal). 
In the case that $M$ has dimension $4$, one key obstruction comes from the signature of $M$, $\sig(M)$, which is an integer-valued invariant attached to $M$ which must vanish whenever $M$ is hyperbolic (see Section \ref{sec:signature} for the definition and details).
One potential strategy to ruling out a hyperbolic structure on an atoroidal surface bundle $M$ is therefore to show that $\sigma(M)\neq 0$.
The main goal of this paper is to show that, on the contrary, there are infinitely many atoroidal surface bundles whose signature vanishes.

We remark that, while all closed hyperbolic manifolds have signature zero, the general relation between signature and negative curvature is a delicate one.
Indeed, Gromov-Thurston \cite{GT87} constructed closed manifolds with curvature pinched close to $-1$ that are not hyperbolic and have signature $0$. For these manifolds, vanishing of the signature can be seen in several different ways. For example, one can use the fact that they are conformally flat and apply Chern-Weil theory, see \cite{Avez}, or one can use Viro's formula for the signature of a branched cover, see \cite{Viro}. Alternatively, a geometric argument can be found in \cite[Thm B]{LafontRoy}. 
On the other hand, Ontaneda \cite[Cor 4]{Ontaneda} obtained closed manifolds with curvature pinched close to $-1$ with non-zero signature.

\vspace{.5cm}

To provide context for the statement of our results, we briefly recall the proof strategy of Theorem \ref{thm:KL}.
In what follows, let $T^2$ denote the $2$-torus and $T^2_X=T^2\setminus X$ for a finite set of points $X$.
We also denote by $M_{4_1}$ the $3$-manifold obtained as the complement of the figure-eight knot in $S^3$. This manifold is well known to admit a unique hyperbolic metric of finite volume and therefore a discrete faithful representation $\pi_1(M_{4_1})\to \PSL_2(\C)=\Isom^+( \Hy^3)$.
The main novelty of \cite{KL} is to show that $\pi_1(M_{4_1})$ admits an injective, type-preserving homomorphism  into $\Mod(T^2_X)$, when $|X|=3$.
By a type-preserving homomorphism, we mean one that takes loxodromic elements to pseudo-Anosov mapping classes, and parabolic elements to reducible mapping classes, where for this classification we have identified $\pi_1(M_{4_1})$ with its image under the aforementioned discrete, faithful representation.
Using branched covers over $T^2$, Kent--Leininger then show that for every $g\ge 4$, there is a finite index subgroup $H_g<\pi_1(M_{4_1})$ with a type-preserving homomorphism of $H_g$ into $\Mod(S_g)$ (see Section \ref{sec:general_goodlift} for more details). 
From there, any closed, quasi-Fuchsian surface subgroup of $H_g$ produces an associated atoroidal surface bundle.

With this strategy in mind, our first result is that the corresponding atoroidal surface bundle has signature $0$, provided the associated closed surface is null-homologous.
In the interest of providing a tool flexible enough to handle any future variations of the construction of Kent--Leininger, we will state and prove it in more generality than just the setting of the previous paragraph.

%%%%%%%%%%%%%%%%%%%%
%%%%%%%%%%%%%%%%%%%%

\begin{thm}\label{thm:general_irregularcover}
Let $S$ be a closed orientable surface, $X\subseteq S$ a finite (possibly empty) set of points, and $S_X=S\setminus X$. Assume that $\chi(S_X)<0$.
Then there exists infinitely many $g\ge 2$ for which there is an injective, type-preserving homomorphism $\iota:\Mod(S_X)\hookrightarrow\Mod(S_g)$.

Moreover, assume that $M$ is a CW-complex  with an injective homomorphism  $\pi_1(M) \hookrightarrow \Mod(S_X)$.
Let $\Sigma\subset M$ be any null-homologous closed orientable surface.
Then the associated $S_g$-bundle $E_\iota$  with base $\Sigma$   has signature $0$.
If $\iota(\pi_1(\Sigma))$ is purely pseudo-Anosov, then $E_\iota$ is atoroidal.
\end{thm}

%%%%%%%%%%%%%%%%%%%%
%%%%%%%%%%%%%%%%%%%%

We remark that the main contribution here is the statement on the vanishing of $\sig(E_\iota)$.
Indeed, the homomorphism $\iota$ in Theorem \ref{thm:general_irregularcover} is constructed via the classical branched covering techniques developed by Birman--Hilden \cite{BH71,BH73}, and more recent developments, see \cite{ALS,KL,MW21}.
In \S\ref{sec:general_goodlift} we propose a variation on this theme that is adapted to our purposes, namely one that gives an embedding of $\Mod(S_X)$ into $\Mod(S_g)$ without having to pass to a finite-index subgroup.

For a concrete application of Theorem~\ref{thm:general_irregularcover}, recall that $H_2(M_{4_1};\Z)=0$ as $M_{4_1}$ is a knot complement in $S^3$ (alternatively a hyperbolic once punctured torus bundle) and therefore every closed, quasi-Fuchsian surface in $M_{4_1}$ is null-homologous.
Applying this theorem to the injective, type preserving homomorphism $\pi_1(M_{4_1})\to\Mod(T^2_X)$ constructed by Kent-Leininger \cite{KL} we conclude the following from Theorem \ref{thm:general_irregularcover}.

%%%%%%%%%%%%%%%%%%%%
%%%%%%%%%%%%%%%%%%%%

\begin{thm}\label{thm:homeotype_2}
For infinitely many natural numbers $h\ge 2$ and infinitely many natural numbers $g\ge 2$, there exists an atoroidal $S_g$-bundle over $S_h$ with signature $0$.
In particular, there are infinitely many homeomorphism types of atoroidal surface bundles over surfaces with signature $0$.
\end{thm}

%%%%%%%%%%%%%%%%%%%%
%%%%%%%%%%%%%%%%%%%%

The observant reader will notice that we only claim infinitely many $g\ge 2$ in contrast to Theorem \ref{thm:KL}.
This is an artifact of our proof and a discussion of this discrepancy can be found in Remark \ref{rem:KLdifference}.

It is known by work of Salter \cite{Salter} that surface bundles over surfaces may have arbitrarily many fibering structures, that is, for every $n$ there is a surface bundle over surface with at least $n$ different fibering structures. 
The constructions of Salter have monodromy contained in the Torelli group, however, Salter \cite{Salter2} has also shown that the stronger property of having monodromy contained in the Johnson kernel implies a unique fibering structure (provided the bundle is non-trivial).
From this perspective, we find the following related question interesting.

%%%%%%%%%%%%%%%%%%%%
%%%%%%%%%%%%%%%%%%%%

\begin{ques}
Given a purely pseudo-Anosov, closed surface subgroup arising from the construction of \cite{KL}, can one describe the monodromy in $\Mod(S_g)$ explicitly? Do the subgroups corresponding to non-commensurable, closed totally geodesic surfaces in $M_{4_1}$ always give rise to distinct homeomorphism types of the associated bundles?
\end{ques}

%%%%%%%%%%%%%%%%%%%%
%%%%%%%%%%%%%%%%%%%%

This would be the case if one could prove that such bundles have a unique fibering structure. \bigskip

%%%%%%%%%%%%%%%%%%%%
%%%%%%%%%%%%%%%%%%%%

\noindent \textbf{Acknowledgements:}
The authors are grateful to Autumn Kent and Chris Leininger for helpful conversations related to their construction of atoroidal surface bundles. Lafont was partially sponsored by NSF DMS--2407438 and Miller was partially sponsored by NSF DMS--2405264.

%%%%%%%%%%%%%%%%%%%%
%%%%%%%%%%%%%%%%%%%%
%%%%%%%%%%%%%%%%%%%%

\section{Background}\label{sec:background}

%%%%%%%%%%%%%%%%%%%%
%%%%%%%%%%%%%%%%%%%%

\subsection{Hyperbolic manifolds and surface bundles}\label{sec:surfacebundle}

%%%%%%%%%%%%%%%%%%%%
%%%%%%%%%%%%%%%%%%%%

We record some well-known facts that we require in the sequel.
In what follows, by a hyperbolic manifold $M$ we always mean a finite-volume quotient of $\Hy^n$ by a discrete, torsion free subgroup of $\Isom^+(\Hy^n)$.
In particular, implicitly all of our manifolds are orientable.

When $n=3$, a particularly important construction of hyperbolic manifolds $M$ are knot complements in the $3$-sphere, which arise from removing a tubular neighborhood of a knot $\mathcal{K}$ from $S^3$ and hyperbolizing the complement.
It is a consequence of Thurston's hyperbolization theorem \cite{Otal,Otal2} that such knot complements are hyperbolizable unless $\mathcal{K}$ is a satellite or torus knot.
As we require it in the sequel, we note that a straightforward Alexander duality computation shows that $\widetilde{H}_k(M;\Z)\cong \widetilde{H}^{2-k}(S^1;\Z)$ and therefore we have the following.

%%%%%%%%%%%%%%%%%%%%

\begin{lem}\label{lem:H2calcknot}
If $M$ is a knot complement in $S^3$, then $H_2(M;\Z)=0$.
\end{lem}

Of particular interest to us in the sequel will be the question of whether certain $4$-manifolds arising as surface bundles over surfaces are hyperbolizable. 
For this, we remind the reader of the classification of $C^\infty$ surface bundles over an arbitrary base $B$, which is due to Earle and Eells \cite{EE} (see also \cite[Prop 4.6]{Morita}).

%%%%%%%%%%%%%%%%%%%%

\begin{thm}[Earle--Eells]\label{thm:ee}
Let $B$ be any $C^\infty$ manifold and $g$ any natural number at least $2$. Then isomorphism classes of $S_g$-bundles over $B$ are in one-to-one correspondence with conjugacy classes of homomorphisms $\pi_1(B)\to\Mod(S_g)$.
\end{thm}

%%%%%%%%%%%%%%%%%%%%%

Indeed, what Earle--Eells show in \cite{EE} is that $\Diff_0(S_g)$ is contractible (provided $g\ge 2$) and therefore the natural map $\Diff(S_g)\to\Mod(S_g)$ is a homotopy equivalence.
This theorem will be useful later in the case that $B$ is either a hyperbolic surface, or a $3$-manifold with boundary.
In the case that $B$ is a closed hyperbolic surface, one can moreover decide when the corresponding closed surface bundle $E$ is atoroidal, i.e. $\pi_1(E)$ has no $\Z^2$ subgroup, using the following result of Bowditch \cite[Lem 1.2]{Bowditch}.

%%%%%%%%%%%%%%%%%%%%%

\begin{lem}\label{lem:atoroidal}
Let $E$ be an $S_g$-bundle over $S_h$ for $g,h\ge 2$. 
Then $E$ is atoroidal if and only if the associated monodromy representation $\pi_1(S_h)\to \Mod(S_g)$ has purely pseudo-Anosov image.
\end{lem}

%%%%%%%%%%%%%%%%%%%%
%%%%%%%%%%%%%%%%%%%%

\subsection{Recollections on signature}\label{sec:signature}

%%%%%%%%%%%%%%%%%%%%
%%%%%%%%%%%%%%%%%%%%

In this subsection, we let $M$ be a closed oriented  manifold of dimension $4k$ for some $k\in\N$.
Then Poincar\'e duality implies that the cup product in middle dimension
$$\cup:H^{2k}(M;\R)\times H^{2k}(M;\R)\to\R,$$
is non-degenerate.
Moreover, as any two classes $\alpha,\beta\in H^{2k}(M;\R)$ have the property that $\alpha\cup\beta=\beta\cup\alpha$, it follows that $-\cup-$ is a non-degenerate, symmetric bilinear form and hence induces a quadratic form $q_M:H^{2k}(M;\R)\to\R$ over $\R$.
By Sylvester's law of inertia, any such quadratic form is isometric to one with $n_+(q_M)$ eigenvalues of $+1$ and $n_-(q_M)$ eigenvalues of $-1$, where $n_+(q_M)+n_-(q_M)$ equals $b_{2k}(M)$, the dimension of $H^{2k}(M;\R)$.
The \emph{signature of $M$} is then given by the formula
$$\sigma(M)=n_+(q_M)-n_-(q_M),$$
which is a homotopy invariant of $M$.
Moreover, we always have the following in the special case that $M$ bounds. 

%%%%%%%%%%%%%%%%%%%%%

\begin{thm}\label{thm:signature}
For a given $4k$-manifold $M$ as above, suppose that there exists a compact oriented  $(4k+1)$-manifold $F$ such that $\partial F=M$, then $\sig(M)=0$.
\end{thm}

%%%%%%%%%%%%%%%%%%%%%

This follows from the usual long exact sequence for the pair $(F,M)$ combined with Poincar\' e duality for manifolds with boundary.
Indeed, in the setting of Theorem \ref{thm:signature} one can show that there is a $(b_{2k}(M)/2)$-dimensional isotropic subspace for $q_M$ from which the result follows.
In fact, one can moreover show that signature is a surjection from the oriented bordism ring $\Omega^{\SO}_{4k}$ to $\Z$ for all $k\in\N$ and is an isomorphism when $k=1$.

We record the following well-known fact about the signature of hyperbolic $4k$-manifolds.
%%%%%%%%%%%%%%%%%%%%

\begin{prop}
If $M$ is a closed hyperbolic $4k$-manifold, then $\sig(M)=0$.
\end{prop}

%%%%%%%%%%%%%%%%%%%%
\begin{proof}
By Hirzebruch's signature formula (see \cite[Thm 8.2.2]{Hirzebruch}), in order to show that the signature vanishes, it  suffices to show that  the Pontryagin classes $p_i(M)$ of $M$ vanish for all $i\in\N$.
For example, in dimension $4$ one has $\sig (M)=\frac 13 <p_1(M),[M]>$.
By the Hirzebruch proportionality principle (see \cite[Appendix 1]{Hirzebruch}), this is equivalent to showing that $p_i(S^{4k-1})=0$.
However this is well known and follows, for instance, from the fact that $TS^{4k-1}$ is stably trivial.
\end{proof}

%%%%%%%%%%%%%%%%%%%%%%%%%%%%%%%%

%%%%%%%%%%%%%%%%%%%%
%%%%%%%%%%%%%%%%%%%%
%%%%%%%%%%%%%%%%%%%%
In the case of $4$-manifolds arising as surface bundles over surfaces, the signature can sometimes be determined via homological data as follows.
Let $\Sigma$ be a closed orientable surface, and assume one is given a homomorphism $\rho:\pi_1(\Sigma)\to \Mod(S_g)$.
This determines an $S_g$-bundle $E_\rho\to\Sigma$, as well as a class $[\rho(\Sigma)]\in H_2(\Mod(S_g);\Z)$, namely the pushforward of the fundamental class of $\Sigma$ along any map $f_\rho:\Sigma\to \operatorname{B}\Mod(S_g)$ inducing $\rho$.

\begin{prop}\label{prop:meyer}
If $g\geq 3$, then
$\sig (E_\rho) =0$ if and only if $[\rho(\Sigma)]=0$.
\end{prop}
\begin{proof}
    This follows from the following facts $$H^2(\Mod(S_g);\Z) \cong \operatorname{Hom}(H_2(\Mod(S_g);\Z),\Z) \cong \Z,$$ 
    $$H_2(\Mod(S_g);\Z) \cong \Z,$$ 
    noting that when $g\geq 3$, $H_1(\Mod(S_g);\Z)=0$, and hence the Ext term vanishes. 
    Also, a generator of $H^2(\Mod(S_g);\Z)$ is given by the Meyer signature cocycle, whose value on $\rho(\Sigma)$ is precisely $\sig (E_\rho)$, see \cite[\S 5.6]{FarbMarg}.
\end{proof}

\begin{rem}\label{rem:pitchforapp}
    Let $M$ be a CW-complex with a non-trivial homomorphism $\rho:\pi_1(M)\to \Mod(S_g)$.
    Suppose that  $M$ admits a map $f:\Sigma\to M$ from a closed orientable surface $\Sigma$ that is $\pi_1$-injective but null-homologous.
    Then, restricting $\rho$ to $\pi_1(\Sigma)$ provides an example of a non-trivial homomorphism $\rho:\pi_1(\Sigma)\to \Mod(S_g)$ for which $[\rho(\Sigma)]=0$.
    The associated surface bundle $E_\rho$ is non-trivial but has trivial signature by Proposition~\ref{prop:meyer}.
    In Appendix \ref{appendix} we provide an alternative proof of this fact by 
    showing that $E_\rho$ is the boundary of a compact orientable smooth $5$-manifold and invoking Theorem~\ref{thm:signature}.    
\end{rem}

%%%%%%%%%%%%%%%%%%%%
%%%%%%%%%%%%%%%%%%%%

\section{Finding type-preserving lifts of \texorpdfstring{$\Mod(S_X)$}{Mod(T2X)}}\label{sec:general_goodlift}
Let $S$ be a closed orientable surface, $X\subseteq S$ a finite set of points, and $S_X=S\setminus X$.
The purpose of this section is to construct injective, type-preserving homomorphisms $\Mod(S_X) \to \Mod(S_g)$ into the mapping class group of a closed surface $S_g$ for some $g\ge 2$.
These can be obtained via the branched covering techniques developed by Birman-Hilden \cite{BH71,BH73,MW21} and more recently in \cite{ALS,IM99,KL}.
In the closed case, i.e., when $X=\varnothing$, this is achieved in \cite[Thm 1, Lem 10]{ALS}, hence we will focus on the punctured case. In the sequel, let $n=|X|\geq 1$.

We start by reviewing how to construct a homomorphism between mapping class groups from branched covers, for the reader's convenience and to fix notation; see \cite[\S2]{ALS} for details.
Let $\kappa:\widetilde{S}\to S$ be a branched cover, branching over $X$.
Let $Y=\kappa^{-1}(X)$ and $\widetilde{S}_Y=\widetilde{S}\setminus Y$, so that $\kappa$ restricts to an unbranched cover $\widetilde{S}_Y\to S_X$ which we also denote by $\kappa$.
If $\Mod_0(S_X)\le \Mod(S_X)$ denotes the finite index subgroup of mapping classes which lift to $\widetilde{S}_Y$ then there is a short exact sequence
%%%%%%%%%%%%%%%%%%%%
\begin{equation}\label{eqn:covering}
\xymatrix{1\ar[r]&K\ar[r]&\widetilde{\Mod}_0(S_X)\ar[r]^-{\kappa_*}&\Mod_0(S_X)\ar[r]&1},
\end{equation}
%%%%%%%%%%%%%%%%%%%%
where $K$ is the deck group of $\widetilde{S}_Y\to S_X$ and $\widetilde{\Mod}_0(S_X)<\Mod(\widetilde{S}_Y)$ is the subgroup of lifts of elements of $\Mod_0(S_X)$.
A common strategy for constructing such a homomorphism is to try to split this sequence, and then compose a section with the map that forgets the punctures $\Mod(\widetilde{S}_Y)\to \Mod(\widetilde{S})$. 
Choices of $\kappa$ for which the sequence \eqref{eqn:covering} splits include the extreme cases of $\kappa$ being characteristic and $\kappa$ having trivial deck group.
In both circumstances, it is possible to obtain the desired embedding into the mapping class group of a closed surface.

\begin{thm}\label{thm:injection_into_closed}
    Let $S$ be a closed orientable surface, $X\subseteq S$ a finite set of points, and $S_X=S\setminus X$.
    Then there is a finite branched cover $\kappa:\widetilde S\to S$, branching over $X$, that induces an injective, type-preserving homomorphism $\Mod(S_X)\to \Mod(\widetilde{S})$. 
    Moreover, the induced unbranched cover 
    $\kappa:\widetilde S_Y = \widetilde{S} \setminus \kappa^{-1}(X) \to S_X$ can be chosen either to be characteristic or to have trivial deck group.
\end{thm}

The main contribution here is that $\Mod(S_X)$ admits an injective, type-preserving homomorphism into the mapping class group of some closed surface $\widetilde{S}$, induced by a cover with a lot of symmetries or no symmetries at all.
Once one has that, one can compose with the injections from \cite[Thm 1]{ALS} to obtain embeddings into higher and higher genera.

A key observation here is that in order for the injections between mapping class groups to be type-preserving, one needs to know that the local degree of $\kappa$ around each branch point is at least $2$; see \cite[Lem 10]{ALS} and \cite[Prop 20]{KL}.

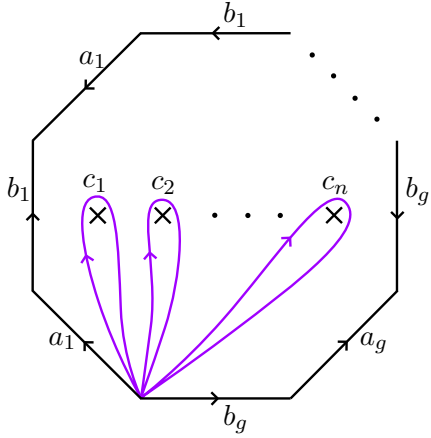
\begin{figure}[t]
\centering
\begin{overpic}[scale=.6]{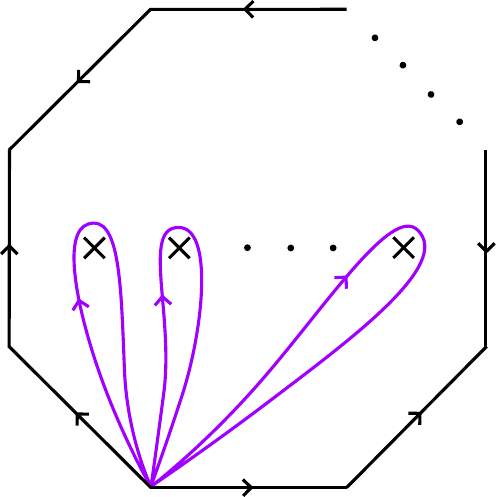}
\put(6,16){$a_1$}
\put(-5,55){$b_1$}
\put(14,90){$a_1$}
\put(52,101){$b_1$}
\put(15,57.5){$c_1$}
\put(33,57){$c_2$}
\put(78,57){$c_n$}
\put(100,55){$b_g$}
\put(88,16){$a_g$}
\put(52,-5){$b_g$}
\end{overpic}
\caption{The fundamental group generators from Equation \eqref{eqn:fundgrp}. }\label{fig:genusg}
\end{figure}

To obtain the desired covers with the required local degree condition, we will study finite covers induced by suitably defined representations of $\pi_1(S_X)$ onto finite groups.
We will use the following presentation for the fundamental group of $S_X$
\begin{equation}\label{eqn:fundgrp}
\pi_1(S_X)=\langle a_1,b_1,\dots,a_g,b_g,c_1,\dots,c_n\mid \prod_{i=1}^g [a_i,b_i]=\prod_{i=1}^n c_i\rangle,
\end{equation}
where the generators are as in Figure \ref{fig:genusg} and $n=|X|\ge 1$ is the number of punctures.
The $c_i$'s will be called \textit{peripheral generators}.
If $g=0$ then it is understood that $a_i,b_i$ do not appear, i.e., all generators are peripheral.

%%%%%%%%%%%%%%%%%%%%
%%%%%%%%%%%%%%%%%%%%
%%%%%%%%%%%%%%%%%%%%
%%%%%%%%%%%%%%%%%%%%

%%%%%%%%%%%%%%%%%%%%
%%%%%%%%%%%%%%%%%%%%

\subsection{Theorem \ref{thm:injection_into_closed} for characteristic covers}
If $\kappa: \widetilde{S}_Y\to S_X$ is a characteristic cover, then the conjugacy class of $\pi_1(\widetilde{S}_Y)$
is invariant under the  action of $\Mod(S_X)$ on  $\pi_1(S_X)$, hence every mapping class lifts, i.e., $\Mod_0(S_X)= \Mod(S_X)$, see \cite[\S 2]{ALS} for details. 

Moreover, one can obtain a splitting for the sequence in Equation \eqref{eqn:covering} following Ivanov-McCarthy \cite[\S 2]{IM99}. 
Indeed, let us choose a point $y_i$ over each branch point $x_i$, for $i=1,\dots, n$.
If $f$ is a homeomorphism of $S_X$, then $f$ permutes the points in $X$ and so we let $x_i=f(x_1)$.
Since the deck group $K$ acts transitively on the fiber over a branch point, up to post-composing with a deck transformation, one can assume that 
$f$ has a lift that sends the chosen $y_1$ to the chosen $y_i$, and such lift is uniquely determined.
This provides a right splitting for the sequence in Equation \eqref{eqn:covering}.

Hence, for any characteristic cover one has an injective homomorphism $\Mod(S_X)\to \Mod(\widetilde{S}_Y)$, that one can then compose with the map that forgets the punctures $\Mod(\widetilde{S}_Y) \to \Mod(\widetilde{S})$.
Therefore to prove Theorem \ref{thm:injection_into_closed} one is reduced to producing a finite characteristic cover with local degree at least $2$ at each puncture.

\begin{lem}\label{lem:characteristic_cover}
Let $\rho:\pi_1(S_X)\to G$ be  a surjective homomorphism to a finite group such that the order of $\rho(c_i)$ in $G$ is at least $2$ for $i=1,\dots,n$.
Then there is a branched cover $\kappa:\widetilde S\to S$ that induces an injective, type-preserving homomorphism $\Mod(S_X)\to \Mod(\widetilde{S})$, and such that
    the associated unbranched cover
    $\kappa:\widetilde S_Y = \widetilde{S} \setminus \kappa^{-1}(X) \to S_X$ is characteristic. 
\end{lem}
\begin{proof}
    The subgroup $\ker(\rho)$ has index $|G|$.
    Since $\pi_1(S_X)$ is finitely generated, there are only finitely many subgroups of index $|G|$, hence intersecting all of them provides a characteristic subgroup $H<\ker(\rho)$ which still has finite index. 
        Let $\kappa:\widetilde{S}_Y\to S_X$ be the unbranched cover corresponding to $H$.
        It can be completed to a branched cover 
        $\kappa:\widetilde{S}\to S$, and from the assumption that the order of $\rho(c_i)$ in $G$ is at least $2$ for $i=1,\dots,n$, the local degree is at least $2$ at each branch point.
    Since $H$ provides a deeper cover, the same is true for $\kappa$.

    By the discussion above, since $\kappa$ is characteristic, we obtain injections
    $\Mod(S_X)\to \Mod(\widetilde{S}_Y) \to \Mod(\widetilde{S)}$, where the last one is the map that forgets the punctures.
\end{proof}

\begin{proof}[Proof of Theorem~\ref{thm:injection_into_closed} for characteristic covers]
By Lemma~\ref{lem:characteristic_cover}, it is enough to construct a representation to a finite group for which the peripheral generators have order at least $2$.
    When $n=|X|\geq 2$, let $\rho:\pi_1(S_X)\to \Z/n\Z$ be the homomorphism $\rho(a_i)=\rho(b_i)=0$ and $\rho(c_i)= 1$. 
When $n=1$, let $\rho:\pi_1(S_X)\to \operatorname{Sym}(3)$ be the homomorphism $\rho(a_i)=\rho(b_i)=0$ for $i\geq 2$, $\rho(a_1)= (12)$, $\rho(b_1)= (23)$ and $\rho(c_1)= [(12),(23)]=(123)$. 
\end{proof}

%%%%%%%%%%%%%%%%%%%%%%%%%%%%%%%%
%%%%%%%%%%%%%%%%%%%%%%%%%%%%%%%%

\subsection{Theorem \ref{thm:injection_into_closed} for totally irregular covers}
In the previous section, we obtained a right splitting in Equation \eqref{eqn:covering} using the fact that the deck group of a characteristic cover acts transitively on each fiber.
On the other end of the spectrum, one can split the sequence in Equation \eqref{eqn:covering} whenever the deck group $K$ of $\kappa$ is trivial.
Indeed, in this case, $\kappa_*$ is an isomorphism. 
The induced injection of mapping class groups is more canonical in the sense that it only depends on the cover $\kappa$ and not a choice of point in the fiber.

To construct such covers from finite groups, the key proposition for us is the following, due to Aramayona--Leininger--Souto \cite[Prop 6]{ALS}, which is stated for closed surfaces but holds equally well for punctured surfaces.
See the discussion after the statement.

%%%%%%%%%%%%%%%%%%%%

\begin{prop}[Aramayona--Leininger--Souto]\label{prop:als}
Suppose that $G$ is a finite group, $\rho:\pi_1(S_X)\to G$ is a surjective homomorphism with characteristic kernel, and $H< G$ is a subgroup such that 
\begin{enumerate}[(a)]
\item $H$ is self-normalizing, i.e., $N_G(H)=H$,
\item $\Aut(G)\cdot H=\Inn(G)\cdot H$.
\end{enumerate}
If $\kappa:\widetilde{S}_Y\to S_X$ is the cover corresponding to $\rho^{-1}(H)<\pi_1(S_X)$. 
Then
\begin{enumerate}[(1)]
\item $\pi_1(\widetilde{S}_Y)$ is invariant under the $\Mod(S_X)$ action on $\pi_1(S_X)$,
\item The sequence in Equation \eqref{eqn:covering} with $\kappa_*$ arising from this cover splits.
\end{enumerate}
In particular, if $H$ satisfies conditions (a), (b), then there is an injective homomorphism $\Mod(S_X)\to \Mod(\widetilde{S}_Y)$. 
\end{prop}

%%%%%%%%%%%%%%%%%%%%

We remark here that condition (a) implies that the deck group $K$ in Equation \eqref{eqn:covering} is trivial (see e.g. \cite[Prop 1.39]{Hatcher}), hence not only does this sequence splits but $\kappa_*$ is actually an isomorphism. 
Condition (b) ensures that every mapping class on $S_X$ lifts to $\widetilde{S}_Y$ and hence $\Mod_0(S_X)$ coincides with $\Mod(S_X)$. 
Indeed, as in the proof of \cite[Prop 1.6]{ALS}, one shows that Condition (b) implies that
$$\Aut(\pi_1(S_X))\cdot \pi_1(\widetilde{S}_Y)=\Inn(\pi_1(S_X))\cdot \pi_1(\widetilde{S}_Y).$$
It follows that $\Out(\pi_1(S_X))$ preserves the conjugacy class of $\pi_1(\widetilde{S}_Y)$ in $ \pi_1(S_X)$. 
Since $\Mod(S_X)\le \Out(\pi_1(S_X))$, it follows as an application of \cite[Prop 1.33]{Hatcher} that every mapping class of $S_X$ lifts to $\widetilde{S}_Y$, i.e., $\Mod_0(S_X)=\Mod(S_X)$.
Therefore in this setting, the  homomorphism in Proposition~\ref{prop:als} is simply induced by the isomorphism $\kappa_*^{-1}:\Mod(S_X)\to  \widetilde{\Mod}_0(S_X)$ composed with the inclusion of the latter into $\Mod(\widetilde{S}_Y)$.

Again, if the local degree of $\kappa$ at each branch point is at least  $2$ then 
the composition $\Mod(S_X)\to \Mod(\widetilde{S}_Y)\to\Mod(\widetilde{S})$ is an injective, type-preserving homomorphism into the mapping class group of a closed surface.

%%%%%%%%%%%%%%%%%%%%

\begin{rem}\label{rem:whatishard}
Constructing examples of pairs $(G,H)$ which satisfy conditions (a), (b) of Proposition \ref{prop:als} is not particularly difficult.
However the additional requirements that one has a surjective homomorphism  $\pi_1(S_X) \to G$ with characteristic kernel and that the associated cover has local degree at least $2$ at each branch point are non-trivial to satisfy.
In the rest of this section we provide pairs $(G,H)$ for which this can be accomplished.
\end{rem}

%%%%%%%%%%%%%%%%%%%%

First, we use similar homomorphisms as in \cite{ALS} to produce examples of covers satisfying the hypotheses of Proposition \ref{prop:als}. 
The following proposition mirrors the argument in the alternate proof of \cite[Thm 1]{ALS}.
However, later on we will have to make different choices for the subgroups and representations involved, so we provide the arguments for the convenience of the reader.

%%%%%%%%%%%%%%%%%%%%

\begin{prop}\label{prop:alspair}
Let $G_0$ be a non-abelian finite simple group, and let $H_0\leq G_0$ be a subgroup satisfying the following assumptions:
\begin{enumerate}
\item $H_0$ is self-normalizing, i.e., $N_{G_0}(H_0)=H_0$.
\item  $\Aut(G_0)\cdot H_0=\Inn(G_0)\cdot H_0$.
\end{enumerate}
Then given any surjective homomorphism $\rho_1:\pi_1(S_X)\to G_0$, there is a $k=k(\rho_1)$, depending on $\rho_1$, such that if we let
$$G=G_0^k=G_0\times\dots\times G_0 \textrm{ and }H=H_0^k=H_0\times\dots\times H_0,$$
then there is a surjective homomorphism $\rho:\pi_1(S_X)\to G$ with characteristic kernel for which the pair $(G,H)$ satisfies the hypotheses of Proposition \ref{prop:als}.
\end{prop}

%%%%%%%%%%%%%%%%%%%%
%%%%%%%%%%%%%%%%%%%%

\begin{proof}
First of all, we show that for any choice of $k$ the pair $(G,H)$ satisfies the hypotheses of Proposition \ref{prop:als}.
Condition (a) follows  from assumption (1).
Moreover, since $G$ is a product of non-abelian finite simple groups, automorphisms of $G$ are compositions of permutations of the factors with automorphisms of each factor. 
Permuting the components of $H$ clearly preserves its conjugacy class, so condition (b) follows from assumption (2).

Now we show that there is a choice of $k$ for which we can construct a surjective homomorphism $\rho:\pi_1(S_X)\to G=G_0^k$ with characteristic kernel.
To this end, let $\rho_1:\pi_1(S_X)\to G_0$ be the given surjective homomorphism.
Considering $\rho_1$ as an element of $\Hom(\pi_1(S_X),G_0)$, we enumerate the $\Aut(\pi_1(S_X))$-orbit of $\rho_1$, acting by precomposition.
$\Aut(G_0)$ also acts on $\Hom(\pi_1(S_X),G_0)$ by postcomposition and we let $\{\rho_1,\dots,\rho_k\}$ be a maximal subcollection of $\Aut(\pi_1(S_X))\cdot \rho_1$ for which no two elements of this subcollection are in the same $\Aut(G_0)$-orbit. We also require that the first element, $\rho_1$, is the same $\rho_1$ we began with.
Then the homomorphism
$$\rho=\rho_1\times\dots\times\rho_k:\pi_1(S_X)\to G=G_0^k,$$
has characteristic kernel by construction and is surjective by a lemma of Hall (see e.g. \cite[Lem 8]{ALS}).
\end{proof}

%%%%%%%%%%%%%%%%%%%%
Now, we enrich the previous proposition by adding conditions on the local degree around the punctures.
Recall that we have fixed the presentation for $\pi_1(S_X)$ in Equation \eqref{eqn:fundgrp}, in which the generators $c_i$ are peripheral loops, each running around a puncture.

\begin{prop}\label{prop:general_goodcovers}
Let $G_0$ be a non-abelian finite simple group, and let $H_0< G_0$ be a subgroup. 
Let $\rho_1:\pi_1(S_X)\to G_0$ be a surjective homomorphism.
For each $i=1,\dots, n$, let $\delta_i$ be the order of $\rho_1(c_i)$ in $G_0$ and let $\delta=\min\{\delta_i\}$.
Assume that:
\begin{enumerate}

\item $H_0$ is self-normalizing.

\item $\Aut(G_0)\cdot H_0=\Inn(G_0)\cdot H_0$.,  

\item For each $i=1,\dots, n$, we have $\delta_i\geq 2$ and  $\gcd( \delta_i,|H_0|)=1$.
\end{enumerate}

Then there exist $k\in \N$ and a branched cover $\kappa:\widetilde{S}\to S$  such that:
\begin{enumerate}[(1)]
\item \label{item:deg} The degree of $\kappa$ is $[G_0:H_0]^k$, and the deck group of $\kappa$ is trivial .

\item \label{item:local_deg} For each $y\in Y=\kappa^{-1}(X)$, the local degree $\deg_y\kappa$ of $\kappa$ at $y$  is at least $\delta$. 
More precisely, if $\kappa (y_i)=x_i\in X$, and
$c_i$ is the peripheral generator corresponding to the puncture $x_i$, then $\delta_i \mid \deg_{y_i}\kappa$.

\item  \label{item:genus}  The genera of $\widetilde{S},S$ satisfy  $g(\widetilde{S})\ge 1+ [G_0:H_0]^k \left(g(S)-1+\dfrac{n(\delta-1)}{2\delta}\right) $.

\item \label{item:injection} The branched cover $\kappa$ induces an injective, type-preserving homomorphism $\iota:\Mod(S_X)\to \Mod(\widetilde{S})$.
\end{enumerate}
\end{prop}

%%%%%%%%%%%%%%%%%%%%
%%%%%%%%%%%%%%%%%%%%

\begin{proof}
Let $k=k(\rho_1)$ be the integer from Proposition \ref{prop:alspair}, $G=G_0^k$, $H=H_0^k$, and let $\rho:\pi_1(S_X)\to G=G_0^k$ be the corresponding surjective homomorphism.
Let $\widetilde{S}_Y\to S_X$ be the unbranched cover corresponding to $\rho^{-1}(H)<\pi_1(S_X)$.
This cover can be completed to a branched cover $\kappa:\widetilde{S}\to S$ branched precisely over points in $X$.
We now proceed to prove all the desired properties of $\kappa$.

To prove \eqref{item:deg}, note that the degree $d$ of $\kappa$ is  the index of $\rho^{-1}(H)$ in $\pi_1(S_X)$, which is in turn equal to the index $[G:H]=[G_0:H_0]^k$. The triviality of the deck group follows from the fact that $N_G(H)=H$.

Let us prove \eqref{item:genus} from \eqref{item:local_deg}.
The Riemann-Hurwitz formula applied to $\kappa$ gives that $\chi(\widetilde S)=d\chi(S)-\sum_{y\in Y} (\operatorname{deg}_y\kappa-1)$.
This says that
\begin{align*}
    g(\widetilde{S})=1+ d(g(S)-1)+\dfrac 12 \sum_{y\in Y} (\operatorname{deg}_y\kappa-1).
\end{align*}
To estimate the correction term, note that the sum of the local degrees over any point  is equal to the degree $d$.
Moreover, \eqref{item:local_deg} implies that there are at most $d/\delta$ points over each $x\in X$. Therefore we have
\begin{align*}
 \sum_{y\in Y} (\operatorname{deg}_y\kappa-1) &=  \sum_{x\in X} \left( \left(\sum_{y\in \kappa^{-1}(x)} \operatorname{deg}_y\kappa\right) - |\kappa^{-1}(x)| \right), \\
 & \ge \sum_{x\in X}\left(d - \dfrac d\delta \right)=\dfrac{dn(\delta-1)}{\delta},
\end{align*}
which proves \eqref{item:genus}.

To prove \eqref{item:injection} from \eqref{item:local_deg}, note that $\widetilde{S}_Y\to S_X$ is by definition the cover corresponding to $\rho^{-1}(H)$. 
By Proposition~\ref{prop:alspair} we have that $\rho$ and $(G, H)$ satisfy the conditions in Proposition~\ref{prop:als}.
Therefore there is an injective homomorphism $\Mod(S_X)\to \Mod(\widetilde{S}_Y)$.
Since by \eqref{item:local_deg} the local degree of $\kappa$ is at least $\delta\geq 2$, we have that this homomorphism is also type-preserving by \cite[Lem 10]{ALS}.
Now consider the homomorphism that forgets the punctures $\Mod(\widetilde{S}_Y) \to \Mod(\widetilde{S})$.
By the second paragraph of the proof of \cite[Prop 20]{KL}, this homomorphism is injective and type-preserving as soon as the local degree of $\kappa$ at each branch point is at least $2$. But once again this is guaranteed by \eqref{item:local_deg}. The homomorphism $\iota$ in \eqref{item:injection} is the composition of these two homomorphisms.

It remains to prove \eqref{item:local_deg}, i.e., to check  the local degree of $\kappa$ at points in $Y$.
Let $y_i\in Y$ and $\kappa (y_i)=x_i\in X$.
Let $c_i$ be the peripheral generator corresponding to the puncture $x_i$.
For any cover of $S_X$, recall that an \textit{elevation} of $c_i$ is any connected component of its full preimage under the covering map. 
Let $\beta_i$ be the elevation of $c_i$ to $\widetilde{S}$ corresponding to $y_i$.
The local degree $\deg_{y_i}\kappa$ is equal to the degree of the cover $\beta_i\to c_i$.
In order to deal with the non-regular cover $\widetilde{S}_Y\to S_X$, we introduce two auxiliary regular covers.
Let $\tau:\widetilde{T}_Z\to S_X$ denote the regular cover corresponding to $\ker(\rho)$. This is a cover with degree $d$ and deck group $G$.
Moreover, there is a cover $\tau':\widetilde{T}_Z\to \widetilde{S}_Y$   such that $\tau= \kappa\circ \tau'$. This cover is also regular, with deck group $H$.

Now, let $\eta_i$ be any elevation of $c_i$ to $\widetilde{T}_Z$. 
Since $\tau:\widetilde{T}_Z\to S_X$ was chosen to be the regular cover corresponding to $\ker(\rho)$, one computes that the degree of the induced cover $\eta_i \to c_i$ is 
\begin{align*}
d_i&=\mathrm{lcm}\{\ord_{\rho_1}(c_i),\ord_{\rho_2}(c_i),\dots,\ord_{\rho_k}(c_i)\},\\
&=\mathrm{lcm}\{\delta_i,\ord_{\rho_2}(c_i),\dots,\ord_{\rho_k}(c_i)\},
\end{align*}
where $\ord_{\rho_i}(c_i)$ is the order of $\rho_i(c_i)$ in $G_0$. 
In particular, note that $\delta_i\mid d_i$.
Let $\beta_i=\tau'(\eta_i)\subseteq \widetilde{S}_Y$. 
Let $d_i'$ be the degree of the induced cover $\eta_i \to \beta_i$. 
Since $\tau'$ is regular with deck group $H$, we have that $d_i'\mid |H|$.

To conclude, note that we have factored the cover $\eta_i\to c_i$ as $\eta_i\to \beta_i\to c_i$. 
Since the degree is multiplicative in covers, we get that $d_i=d'_i\deg_{y_i}\kappa$. 
As noted, $\delta_i\mid d_i$ and $d_i'\mid |H|$.
It follows that $\delta_i \mid |H|\deg_{y_i}\kappa=|H_0|^k\deg_{y_i}\kappa$.
Finally, the assumption that $\delta_i$ and $|H_0|$ are coprime implies that $\delta_i\mid \deg_{y_i}\kappa$, as desired.
\end{proof}

%%%%%%%%%%%%%%%%%%%%%%%%%%%%%%%%%%%%%%%%%%
%%%%%%%%%%%%%%%%%%%%%%%%%%%%%%%%%%%%%%%%%%

The remainder of this section will be devoted to producing covers satisfying Proposition \ref{prop:general_goodcovers}, which is the content of the following statement, and will complete the proof of Theorem \ref{thm:injection_into_closed} for covers with trivial deck group. 

\begin{thm}\label{thm:goodlifts}
    Let $S$ be a closed orientable surface, $X\subseteq S$ a finite set of $n$ points, and $S_X=S\setminus X$. Assume $\chi(S_X)<0$.
    Let $p\geq \max\{n,13\}$ be a prime such that $p\equiv 1 \pmod{4}$.
    Then there exists $k\in \N$ and a branched cover $\kappa:\widetilde{S}\to S$ such that
    \begin{enumerate}
        \item The branched cover $\kappa$ induces an injective, type-preserving homomorphism $\iota:\Mod(S_X)\to \Mod(\widetilde{S})$.
        
        \item The genus of $\widetilde{S}$ satisfies 
        $g(\widetilde{S})\ge 1+ (p+1)^k \left(g(S)-1+\dfrac{n}{2}\dfrac{p-1}{p+1}\right)$.

        \item The deck group of $\kappa$ is trivial.
    \end{enumerate}
\end{thm}
\begin{proof}
    This boils down to producing a pair $(G_0,H_0)$ and a surjective homomorphism $\rho_1:\pi_1(S_X)\to G_0$ satisfying the assumptions of Proposition~\ref{prop:general_goodcovers}.
    
    To accomplish this, we need to split into cases depending on the topology of $S_X$.
    We refer the reader to the next three subsections \S\ref{sec:generic}, \ref{sec:oncepunctured}, \ref{sec:spheres} for explicit constructions, where $G_0$ is always chosen to be $G_0=\PSL_2(\F_p)$.
    On the other hand, we change $H_0$ and $\rho_1$. 
    Our choice of $H_0$ will always have order at most $p(p-1)/2$, which gives the smallest possible value of the degree.
    Moreover, the order of a peripheral generator will always be at least $(p+1)/2$, which gives the values in the genus computation.
\end{proof}

%%%%%%%%%%%%%%%%%%%%%%%%%%%%%%
%%%%%%%%%%%%%%%%%%%%%%%%%%%%%%

\subsubsection{The generic case}\label{sec:generic}
We start with the case $g\geq 1$ and $n\geq 2$. 
We use the  presentation for the fundamental group of $S_X$ in Equation~\eqref{eqn:fundgrp},
where the generators are as in Figure \ref{fig:genusg}.

For the convenience of the reader, we note that the group $G_0=\PSL_2(\F_p)$ is simple for any prime $p\ge 5$ and that $\SL_2(\Z)$, and hence its quotients, is generated by the two matrices 
$\begin{pmatrix}
1&1\\
0&1\end{pmatrix}$
and 
$\begin{pmatrix}
1&0\\
1&1\end{pmatrix}$.

Let $A_0<G_0$ the subgroup of diagonal matrices, and $H_0$ the normalizer of $A_0$ given explicitly as
$$H_0=N_{G_0}(A_0)=\left\{h\in G_0~\middle|~ h=\begin{pmatrix}
a&0\\
0&1/a\end{pmatrix}\text{ or }\begin{pmatrix}
0&-a\\
1/a&0\end{pmatrix}\text{ with }a\in\mathbb{F}_p^*\right\}.$$

To verify condition (1) in Proposition~\ref{prop:alspair}, note that $H_0$ is a maximal subgroup of $G_0$ (\cite[Prob III.\S6.7]{Suzuki}) and hence is self-normalizing in $G_0$ as $G_0$ is a simple group.
To verify condition (2), note that $\Aut(G_0)\cong\PGL_2(\mathbb{F}_p)$ and moreover any automorphism of $G_0$ is a composition of an inner automorphism and conjugation by the diagonal matrix 
$$d_0=\begin{pmatrix}
1&0\\
0&\epsilon
\end{pmatrix},$$
where $\epsilon\in\F_p^*\setminus \F_p^{*2}$ is any fixed non-square number (see for instance \cite[Thm 30]{Steinberg}).
One computes that conjugation by $d_0$ normalizes $H_0$ and therefore $\Aut(G_0)\cdot H_0=\Inn(G_0)\cdot H_0$. 

Define a representation $\rho_1:\pi_1(S_X)\to G_0$ as follows:
\begin{align*}
    \rho_1(a_1)=\rho_1(b_1)=\begin{pmatrix}
    1&1\\
    0&1
    \end{pmatrix}, & \quad \rho_1(a_i)=\rho_1(b_i)=\begin{pmatrix}
    1&0\\
    0&1
    \end{pmatrix}, \textrm{ for } i\geq 2\\
    \rho_1(c_i)=\begin{pmatrix}
    1&0\\
    1&1
    \end{pmatrix} \textrm{ for } i< n, & \quad \rho_1(c_n)=\begin{pmatrix}
    1&0\\
    p-n+1&1
    \end{pmatrix}.
    \end{align*}
A direct computation shows that for any prime $p\ge n$  the above assignments give a well defined surjective homomorphism $\rho_1:\pi_1(S_X)\to G_0$ such that $\rho_1(c_i)$ has order $p$ in $G_0$ for all $i=1,\dots,n$. 
In particular, note that the order of $\rho_1(c_i)$ is coprime with $|H_0|=p-1$.

%%%%%%%%%%%%%%%%%%%%%%%%%%%%%%%%%%%%%%%%%%
%%%%%%%%%%%%%%%%%%%%%%%%%%%%%%%%%%%%%%%%%%

\subsubsection{The once-punctured case}\label{sec:oncepunctured}
 If  $n=1$ and $g\ge 1$, then we need to change the representation $\rho_1:\pi_1(S_X)\to G_0$ used in the case $n\geq 2$. 
Let $C$ be any element of order $(p+1)/2$ in $G_0$.
 Such elements exist because $G_0$ contains a dihedral group of order $p+1$.
Given any $p\ge 13$, the main theorem in \cite{MW11} guarantees the existence of two elements $A,B\in G_0$  such that $\langle A,B\rangle=G_0$ and $[A,B]=C$.
    Then, still referring to the presentation in Equation~\eqref{eqn:fundgrp}, we set 
    \begin{align*}
    &\rho_1(a_1)=A, \quad \rho_1(b_1)=B, \quad \rho_1(c_1)=C\\
    &\rho_1(a_i)=\rho_1(b_i)=\begin{pmatrix}
    1&0\\
    0&1
    \end{pmatrix}, \textrm{ for } i\geq 2.
    \end{align*}
    This gives a well defined surjective homomorphism $\rho_1:\pi_1(S_X)\to G_0$ such that $\rho_1(c_1)$ has order $(p+1)/2$.

    Let $H_0$ be the subgroup of upper triangular matrices.  
    As shown in \cite{ALS}, the pair $(G_0,H_0)$ satisfies the assumptions (1) and (2) in  Proposition~\ref{prop:general_goodcovers}.
    Moreover, $H_0$ has order $p(p-1)/2$, which is coprime with $(p+1)/2$, as $(p+1)/2$ and $(p-1)/2$ are consecutive integers.

%%%%%%%%%%%%%%%

\subsubsection{The genus zero case}\label{sec:spheres}
Let $S_{0,n}$ denote the sphere with $n$ punctures and present the fundamental group as in Equation \eqref{eqn:fundgrp} so that
$$\pi_1(S_{0,n})=\langle c_1,\dots, c_n\mid \prod_{i=1}^nc_i=1\rangle.$$
Suppose that $p>n-2$, so that $(n-2)$ has a multiplicative inverse $s$ in $\F_p$.
Moreover, let $t\in\F_p$ be an element so that any solution to the polynomial $x^2-(2+t)x+1$ in $\F_{p^2}$ has order $p+1$, necessarily some such $t\in\F_p$ exists.
Define a representation $\rho_1:\pi_1(S_{0,n})\to G_0$ by
\begin{align*}\rho_1(c_i)&=\begin{pmatrix}
1&0\\
s&1
\end{pmatrix},\text{ for }i\le n-2,\quad\quad
\rho_1(c_{n-1})=\begin{pmatrix}
1&t\\
0&1\end{pmatrix},\\
\rho_1(c_{n})&=\left(\begin{pmatrix}
1&0\\
1&1\end{pmatrix}\begin{pmatrix}
1&t\\
0&1\end{pmatrix}\right)^{-1}=\begin{pmatrix}
1+t&-t\\
-1&1
\end{pmatrix}
\end{align*}
Note that $\rho_1$ is surjective since if $k$ is the multiplicative inverse of $t$ in $\F_p$, then
$$\rho_1(c_1\dots c_{n-2})=\begin{pmatrix}
1&0\\
1&1\end{pmatrix},\quad\quad \rho_1(c_{n-1})^k=\begin{pmatrix}
1&1\\
0&1\end{pmatrix},$$
which generate $G_0$.
Moreover by construction $\rho_1(c_i)$ has order $p$ for all $1\le i\le n-1$. 
Note that $\mathrm{trace}(\rho_1(c_n))=2+t$ and therefore this element is non-parabolic since $t\neq 0$ \cite[Lem 3.2]{MW11}.
Moreover, by our choice of $t$, $\rho_1(c_n)$ has order $(p+1)/2$ by examining its characteristic polynomial and noting that we are in $\PSL_2(\F_p)$ not $\SL_2(\F_p)$.

Letting $H_0$ be as in the generic case of \S\ref{sec:generic}, we have that $|H_0|=p-1$, which in particular is coprime to both $p$ and $(p+1)/2$ whenever $p\equiv 1\pmod{4}$.

\begin{rem}
The attentive reader will notice that such a representation actually defines a branched cover for any $S_{g,n}$ provided $n\ge 3$ by considering the composite map
$$\pi_1(S_{g,n})\twoheadrightarrow\pi_1(S_{0,n})\to G_0,$$
where the former map sends all $a_i$, $b_i$ to the identity element. In particular, this gives another way of using Proposition~\ref{prop:general_goodcovers} in that special case.

We also point out that we did not actually need to choose $t$ satisfying the above condition, provided one is willing to choose $H_0$ based on $t$.
More explicitly, one can choose any $t\in\F_p^*$ and define the similar representation as above.
If $x^2-(2+t)x+1$ is reducible then $\rho_1(c_n)$ will have order bigger than $1$, due to the trace condition, and dividing $(p-1)/2$. In that case, one can choose $H_0$ to be the dihedral group of order $p+1$ and still conclude since the orders of $\rho_1(c_i)$ and $|H_0|$ are coprime.
If $x^2-(2+t)x+1$ is irreducible then $\rho_1(c_n)$ will have order bigger than $1$ and dividing $(p+1)/2$ and one similarly chooses $H_0$ to be the dihedral group of order $p-1$ and concludes as above.
In particular, any non-zero choice of $t$ gives an appropriate branched cover.
\end{rem}

\section{Proofs of Theorems \ref{thm:general_irregularcover} and \ref{thm:homeotype_2}}

In this section, we combine the ingredients from Sections \ref{sec:background} and \ref{sec:general_goodlift} to finish the proofs of the main theorems.

%%%%%%%%%%%%%%%%%%%%
%%%%%%%%%%%%%%%%%%%%
%%%%%%%%%%%%%%%%%%%%

\begin{proof}[Proof of Theorem \ref{thm:general_irregularcover}]
From Theorem~\ref{thm:injection_into_closed}, there is an injective, type-preserving homomorphism
$$\iota: \Mod(S_X)\to \Mod(\widetilde{S}),$$
into the mapping class group of a closed orientable surface $\widetilde{S}$.
Composing with the injections from \cite[Thm 1]{ALS}, one obtains injective, type-preserving homomorphisms $\iota:\pi_1(M)\hookrightarrow\Mod(S_g)$ for arbitrarily high genus $g$.
Alternatively, one can use Theorem~\ref{thm:goodlifts} for an increasing sequence of primes.

Restricting $\iota$ to $\pi_1(M)$, and then to $\pi_1(\Sigma)$, we have produced injective, type-preserving homomorphisms $\iota:\pi_1(\Sigma)\to \Mod(S_g)$ for infinitely many values of $g$.
Corresponding to each of these homomorphisms, Theorem \ref{thm:ee} provides an atoroidal $S_g$-bundle $E_\iota$ over $\Sigma$.
Since  $[\Sigma]=0$ in $H_2(M;\Z)$, we also have that $[\iota_p(\Sigma)]=0$ in $H_2(\Mod(S_g);\Z)$ and therefore we conclude that $\sig(E_\iota)=0$ by Proposition~\ref{prop:meyer}.
Alternatively, one can use the results in Appendix \ref{appendix} to prove directly that $E_\iota$ bounds a compact oriented $5$-manifold, and then invoke Theorem~\ref{thm:signature}.

The final statement is a consequence of
 Lemma \ref{lem:atoroidal} and the fact that $\iota$ is type-preserving.
\end{proof}

%%%%%%%%%%%%%%%%%%%%
\begin{proof}[Proof of Theorem \ref{thm:homeotype_2}]
    By work of Kent-Leininger \cite{KL}, 
    the fundamental group of the figure-eight knot complement $M_{4_1}$ admits an injective, type-preserving representation into $\Mod(T^2_X)$, where $T^2$ is the torus and $|X|=3$.

Let $\Sigma \cong S_h$ be any closed-orientable quasi-Fuchsian surface in $M_{4_1}$ of genus $h$. 
In particular, $\Sigma$ is $\pi_1$-injective and totally loxodromic; moreover, $\Sigma$ is null-homologous by Lemma~\ref{lem:H2calcknot}.
    By Theorem~\ref{thm:general_irregularcover}
    we obtain infinitely many $g$ for which there exist an atoroidal $S_g$-bundle $E_{g,h}$ with base $\Sigma \cong S_h$ with signature $0$.

Note that $\chi(E_{g,h})=4(g-1)(h-1)$, so for fixed $h$ and different values of $g$ we obtain non-homeomorphic bundles.   \end{proof}

\begin{rem}\label{rem:KLdifference}
In the case of $S=T^2$ and $|X|=3$, Kent-Leininger \cite{KL} construct an injective, type-preserving homomorphism from a finite index subgroup of $\Mod(T^2_X)$ into $\Mod(S_g)$ for all $g\geq 4$.
This then provides a homomorphism $\iota$ as above for a finite cover of the figure-eight knot complement.
As a null-homologous surface does not necessarily lift to a null-homologous surface, for the purposes of Theorem~\ref{thm:general_irregularcover} we do not want to pass to finite index subgroups.
So here we used a variation of the techniques from \cite{ALS,KL} that work only for infinitely many $g$, but do not require to pass to a finite index subgroup.
\end{rem}

%%%%%%%%%%%%%%%%%%%%
%%%%%%%%%%%%%%%%%%%%
%%%%%%%%%%%%%%%%%%%%

\appendix
\section{Appendix: bounding null-homologous surfaces}\label{appendix}
The purpose of the appendix is to provide a proof of the following proposition that does not rely on Proposition~\ref{prop:meyer}, as promised in  Remark~\ref{rem:pitchforapp}.

\begin{prop}\label{prop:without_meyer}
    Let $M$ be a CW-complex, $\rho:\pi_1(M)\to \Mod(S_g)$ a homomorphism, and let $f:\Sigma\to M$ be a map from a closed orientable surface that is $\pi_1$-injective but null-homologous. 
    Let $E$ be the $S_g$-bundle over $\Sigma$ associated with the restriction of $\rho$ to $f_*(\pi_1(\Sigma))$.
    Then there exists a compact orientable smooth $5$-manifold $W$ such that $\partial W=E$. In particular, $\sig(E)=0$.
\end{prop}

This is based on the fact that, while in general
null-homologous classes need not bound genuine manifolds,  in low dimensions this is the case. 
More precisely, we have the following statement.

%%%%%%%%%%%%%%%%%%%%

\begin{lem}\label{lem:bordism}
Let $f: \Sigma \rightarrow M$ be a continuous map from a closed
oriented surface into a CW-complex $M$. If $f_*[\Sigma]=0 \in H_2(M; \mathbb Z)$, then $\Sigma$ bounds a connected
oriented $3$-manifold $W$, and there exists a continuous map $F:W\rightarrow M$ with the property that $F|_{\partial W} = f$.
\end{lem}

%%%%%%%%%%%%%%%%%%%%

While this result is well-known, we outline two different proofs for the convenience of the reader. The first proof uses a high level argument on bordism theory, while the second one is more elementary and provides a concrete construction.

\begin{proof}[First proof of Lemma \ref{lem:bordism}]
We recall Atiyah's bordism theory for a topological space $M$ (see \cite{Atiyah}). 
Elements of the bordism group $\MSO_k(M)$ consist of equivalence classes of pairs $(N, f)$, where 
$N$ is a closed (but not necessarily connected) oriented smooth $k$-manifold, and $f$ is a continuous map
$f: N \rightarrow M$.
Two elements $(N_1, f_1)$ and $(N_2, f_2)$ are equivalent if there exists a compact oriented smooth $(k+1)$-manifold $W^{k+1}$ satisfying 
$\partial W^{k+1} = N_1 \coprod (-N_2)$ and a continuous map $F: W \rightarrow M$ satisfying $F|_{N_i}= f_i$. The group operation is given by disjoint union,
and the zero element is the pair $(\emptyset, \emptyset)$. In particular, the lemma is equivalent to the statement that $(\Sigma, f) = 0 \in \MSO_2(M)$. 

Note that the functor $\MSO_*$ forms a generalized homology theory \cite{Atiyah} -- it satisfies all the axioms of homology, except for the dimension axiom. Indeed, the $\MSO$-groups of a point are the oriented bordism groups $\Omega^{\SO}_q$, which in low dimensions satisfy $\Omega^{\SO}_0\cong \mathbb Z$, $\Omega^{\SO}_q=0$ for $1\leq q\leq 3$, and $\Omega^{\SO}_4\cong \mathbb Z$ with the isomorphism provided by the signature.
Moreover there is a natural homomorphism $\MSO_k(M) \rightarrow H_k(M; \mathbb Z)$, given by $(N, f) \mapsto f_*[N]$.

To complete the proof, we just verify that the natural homomorphism is an isomorphism in degrees $k\leq 3$. 
For the CW-complex $M$, one
can compute the generalized homology group $\MSO_k(M)$ by using an Atiyah--Hirzebruch type spectral sequence (see for instance \cite[Thm 9.6]{DK}). 
The $E^2$-term is given by
$$E^2_{p,q} := H_p(M; \Omega^{\SO}_q),$$
and the spectral sequence converges to $\MSO_{p+q}(M)$. 

Since the coefficients $\Omega^{\SO}_q$ vanish for $q<0$ and homology vanishes for $p<0$, this is a first quadrant spectral sequence. 
The first few rows are:
$$
E^2_{p,q} = \begin{cases} 
H_p(M; \mathbb Z), & q=0, 4 \\
0, & 1\leq q \leq 3. \\
\end{cases}
$$
In particular, the only non-zero entries satisfying $p+q\leq 3$ are the four entries $E^2_{p,0} = H_p(M; \mathbb Z)$ for $0\leq p \leq 3$.
The differentials in the spectral sequence take the form $d^r_{p,q} : E^r_{p,q} \rightarrow E^r_{p-r, q+r-1}$.
This immediately implies that, within the
range $p+q\leq 3$, we have $E^\infty_{p, q} = E^2_{p,q}$. Putting this all together, we obtain the equalities:
$$\MSO_p(M) = E^\infty _{p,0} = E^2_{p,0} = H_p(M;\mathbb Z),$$
within the range $p\leq 3$. Since the natural homomorphism $\MSO_p(M)\rightarrow H_p(M;\mathbb Z)$ coincides with the projection to the $E^\infty_{p,0}$ term, 
this completes the proof of the lemma.
\end{proof}

%%%%%%%%%%%%%%%%%%%%

Alternatively, one can give a constructive proof of Lemma \ref{lem:bordism}, which we now sketch.
In the following, an $n$-dimensional pseudomanifold with boundary is a simplicial complex $Y$ such that every simplex of $Y$ is contained in at least one $n$-simplex, and every ($n-1$)-simplex of $Y$ is contained in exactly one or two $n$-simplices. The subcomplex consisting of ($n-1$)-simplices that are contained in exactly one $n$-simplex is called the boundary of $Y$, and it is denoted by $\partial Y$. 
When $\partial Y=\varnothing$, we say that $Y$ is a closed pseudomanifold. 

%%%%%%%%%%%%%%%%%%%%

\begin{proof}[Second proof of Lemma \ref{lem:bordism}]
Take a triangulation of $\Sigma$, and consider the $2$-chain $\sum (f\circ \sigma_i) \in C_2(M; \mathbb Z)$, where $\sigma_i$ are the triangles in the triangulation
of $\Sigma$. Since this chain represents $f_*[\Sigma]=0$, there is a $3$-chain $\sum b_i \tau_i \in C_3(M; \mathbb Z)$ satisfying
$$\partial \big( \sum b_i \tau_i  \big) = \sum (f\circ \sigma_i).$$
After potentially composing with some orientation reversing involution of the standard simplex $\Delta^3$, we may assume all the coefficients $b_i$ are positive. Since all singular $2$-simplices
in $\sum b_i \partial (\tau_i)$ cancel, with the exception of those in $\sum (f\circ \sigma_i)$, we can form an abstract CW-complex by taking $\sum b_i$ tetrahedra,
and pairwise identifying facets that cancel out. After doing these identifications, one obtains a $3$-dimensional compact oriented pseudomanifold with boundary $Y$ and a map 
$F:Y\rightarrow M$ induced by the individual maps $\tau_i:\Delta ^3 \rightarrow M$ on the individual tetrahedra (see \cite[Pgs 108--109]{Hatcher}, 
or \cite[Prop 2.1]{LP}). Note that by construction, $\partial Y$ can be 
identified with the triangulated surface $\Sigma$, and the map $F$ restricts to $f$ on $\partial Y = \Sigma$. 

Without loss of generality, we can assume $Y$ is connected. The links of lower dimensional facets in $Y$ are themselves connected oriented pseudomanifolds with boundary of the appropriate dimension. 
Now consider the possible non-manifold points in $Y$. 
Along any $2$-face, we have exactly one or two tetrahedra glued together. The points in the interior of a $2$-face are respectively manifold boundary points or manifold points.
Along any edge, the link in $Y$ will be a connected $1$-dimensional pseudomanifold with boundary. 
There are only two possibilities: an interval (if the edge lies in $\partial Y$) or a circle (if the edge is not on $\partial Y$). The points in the interior of an edge are respectively manifold boundary points or manifold points.
Finally, we can consider the links of vertices. Again, there are two possibilities. 
If $v$ is a vertex not on $\partial Y$, then the link is a connected oriented closed $2$-dimensional pseudomanifold, hence an oriented surface $\Sigma_g$ of some genus $g\geq 0$. The vertex $v$ will be a manifold point if and only if the link is a sphere $S^2 = \Sigma_0$. 
Similarly, if $v$ is a vertex on $\partial Y$, then the link will be a connected oriented $2$-dimensional pseudomanifold with boundary, with a single $S^1$ boundary component (the link of $v$ in $\partial Y$), i.e. a surface $\Sigma _ {g,1}$ of genus $g$ with a single boundary component. Again, $v$ will be a manifold boundary point if and only if the link is a disk $\mathbb D^2 = \Sigma_{0,1}$. 

Finally, we resolve the singularities of $Y$ by forming a $3$-manifold $\hat Y \rightarrow Y$. $\hat Y$ is formed by modifying a small neighborhood of each 
non-manifold vertex of $Y$. First consider vertices $v$ that do not lie on $\partial Y$. If the link is $\Sigma _g$ of genus $g\geq 1$, you replace the 
conical neighborhood $C(\Sigma_g)$ 
by gluing in a handlebody $H_g$ of genus $g$ along $\partial H_g = \Sigma_g = \partial C(\Sigma_g)$. Collapsing a tubular neighborhood of $\partial H_g$ gives
an induced map $H_g \rightarrow C(\Sigma_g)$. Next, consider vertices that lie on $\partial Y$, with link some surface $\Sigma_{g,1}$ with boundary component
a circle. Again, one can replace the conical neighborhood $C(\Sigma_{g,1})$ by gluing in a handlebody $H_g$, where now you attach the handlebody by identifying 
a smoothly embedded $\Sigma_{g,1} \subset \Sigma_g = \partial H_g$ with the subset $\Sigma_{g,1} \subset C(\Sigma_{g,1})$. The subset 
$\Sigma_{g,1} \subset H_g$ we are gluing along has a small tubular neighborhood in $H_g$ which is homeomorphic to $\Sigma_{g,1} \times [0, 1]$, so collapsing 
the complement of this neighborhood gives a map $H_g \rightarrow C(\Sigma_{g,1})$. After making all these replacements, we obtain the connected 
oriented $3$-manifold with boundary $\hat Y$ and the (neighborhood collapsing) map $\psi: \hat Y \rightarrow Y$. By construction, the map $\psi$ has degree one, i.e. $\psi_*[\hat Y] = [Y]$, and is homotopic to a homeomorphism on the boundary. Finally, we get the desired bounding map $\hat F = F \circ \psi : \hat Y \rightarrow Y \rightarrow M$.
\end{proof}

%%%%%%%%%%%%%%%%%%%%

\begin{rem}
We note that Lemma \ref{lem:bordism} can fail if we replace $\Sigma$ by a closed oriented manifold $N$ of dimension $\geq 4$. For a concrete example, consider the cone $M= C(\mathbb C P^2)$, $N=\mathbb CP^2$, and the map $f: N \rightarrow M$ given by the inclusion of $N$ as the base of the cone. 
Since $M$ is contractible, we have $H_4(M;\Z)=0$, hence clearly $f_*[N]=0$. However, since the signature of $N=\mathbb CP^2$ is non-zero, there are no compact orientable $5$-manifolds $W$ that are bounded by $N=\mathbb CP^2$.  

In the spectral sequence proof, the issue is that in degree $4$, the spectral sequence gives rise to a short exact sequence
$$\xymatrix{
0\ar[r]& E^\infty_{0,4} \ar[r]& \MSO_4(M) \ar[r]& E^\infty_{4,0} \ar[r]& 0}.$$
The term $E^\infty_{4,0} = E^2_{4,0} = H_4(M; \Z)$ shows that there is still a surjection to $H_4(M; \Z)$, but it now has a kernel $E^\infty_{0,4} = E^2_{0,4} =
H_0(M; \Omega^{\SO}_4) \cong \Omega^{\SO}_4 \cong \Z$. The example described above provides a concrete element in the kernel.

In the constructive proof, one can start with a triangulation of the $4$-manifold $N$ satisfying $[N]=0 \in H_4(M;\Z)$ and construct a $5$-dimensional pseudomanifold $Y$ with boundary $\partial Y =N$. However when one tries to resolve the singularities of $Y$ to turn it into a manifold $W=\hat Y$, we eventually end up having
to resolve singularities at vertices. At an interior vertex $v$, the link is a closed oriented $4$-manifold $L_v$. If the signature of $L_v$ is zero, one can
remove the conical neighborhood $C(L_v)$ and replace it by a bounding manifold $M_v$ (with $\partial M_v = L_v$). In the general case, we might not
be able to resolve all of these singularities. The obstruction to doing this lies in the sum of the signatures of the $L_v$, which one can naturally view as an 
element in $H_0(Y; \Omega^{\SO}_4)\cong \mathbb Z$. Thus the failure of the constructive approach naturally gives rise to the same obstruction 
class as the spectral sequence approach.
\end{rem}

%%%%%%%%%%%%%%%%%%%%

As an immediate corollary, we obtain the following.

%%%%%%%%%%%%%%%%%%%%

\begin{cor}\label{cor:bounding}
Suppose $M$ is a manifold and $\Sigma$ is a closed, essential surface in $M$ such that $[\Sigma]=0\in H_2(M;\Z)$.
Then there is a connected oriented $3$-manifold $N$ with $\partial N=\Sigma$ and a continuous map $F:N\to M$ such that $F\vert_{\partial N}=\id_\Sigma$.
\end{cor}

We can now prove the main proposition of this appendix. 

\begin{proof}[Proof of Proposition~\ref{prop:without_meyer}]
Since $f:\Sigma \to M$ is null-homologous, 
Corollary \ref{cor:bounding} provides a connected oriented $3$-manifold $N$ with $\partial N=\Sigma$ and a continuous map $F:N\to M$ such that $F\vert_{\partial N}=\id_{\Sigma}$.
In particular, $\partial N=\Sigma$ is $\pi_1$-injective in $N$, and
we have the inclusions
$$\pi_1(\Sigma)= F_*\left(\pi_1(\partial N)\right)\le F_*(\pi_1(N))<\pi_1(M).$$
Consider the composite homomorphism
$$\iota\circ F_*:\pi_1(N)\to \Mod(S_g).$$
Then  by Theorem \ref{thm:ee} there is an $S_g$-bundle $W$ over $N$ whose boundary $\partial W$ is an $S_g$-bundle over $\partial N=\Sigma$. 
Since $F\vert_{\partial N}=\id_\Sigma$, the monodromy restricted to $\partial N$ is identical to that of $E$ and therefore Theorem \ref{thm:ee} shows that $\partial W$ and $E$ are isomorphic as $S_g$-bundles over $S_h$. In particular, $\partial W$ and $E$ are diffeomorphic.
By Theorem \ref{thm:signature}, $\sig(\partial W)=0$ hence we conclude that $\sig(E)=0$ as well, as required.
\end{proof}

%%%%%%%%%%%%%%%%%%%%
%%%%%%%%%%%%%%%%%%%%
%%%%%%%%%%%%%%%%%%%%

%%%%%%%%%%%%%%%%%%%%
%%%%%%%%%%%%%%%%%%%%
%%%%%%%%%%%%%%%%%%%%
\bibliography{signature}
%%%%%%%%%%%%%%%%%%%%
\Addresses

%%%%%%%%%%%%%%%%%%%%
\end{document}